# A Data-Driven, Energy-based Approach for Identifying Equations of Motion in Vibrating Structures Directly from Measurements

Cristian López[1], Aryan Singh[1], Ángel Naranjo[1], Keegan J. Moore[1,2]*

*[1]Department of Mechanical and Materials Engineering, University of Nebraska-Lincoln, Lincoln, NE 68588*
*[2]Daniel Guggenheim School of Aerospace Engineering, Georgia Institute of Technology, Atlanta, GA 30332*

*Corresponding author: (K.J. Moore)
E-mail address: kmoore@gatech.edu

## ABSTRACT

Determining the underlying equations of motion and parameter values for vibrating structures is of great concern in science and engineering. This work introduces a new data-driven approach called the energy-based dual-phase dynamics identification (EDDI) method for identifying the nonlinear dynamics of single-degree-of-freedom oscillators. The EDDI method leverages the energies of the system to identify the governing dynamics through the forces acting on the oscillator. The approach consists of two phases: a model-dissipative and model-stiffness identification. In the first phase, the fact that kinetic and mechanical energies are equivalent when the displacement is zero is leveraged to compute the energy dissipated and a corresponding model for the nonlinear damping of the system. In the second phase, the energy dissipated is used to compute the mechanical energy (ME), which is then used to obtain a reformulated Lagrangian. The conservative forces acting on the oscillator are then computed by taking the derivative the Lagrangian, then a model for the nonlinear stiffness is identified by solving a system of linear equations. The resulting governing equations are identified by including both the nonlinear damping and stiffness terms. A key novelty of the EDDI method is that the only thing required to perform the identification is free-response measurements and the mass of the oscillator. No prior understanding of the dynamics of the system is necessary to identify the underlying dynamics, such that the EDDI method is a truly data-driven method. The method is demonstrated using simulated and measured responses of nonlinear single-degree-of-freedom systems with a variety of nonlinear mechanisms.

## 1. Introduction

Vibrating structures hold a pivotal role across engineering and scientific domains, influencing the design and functionality of many systems and devices [1,2]. An understanding of the governing dynamics of these structures is paramount to optimize their performance and guarantee their reliability and safety [3]. Historically, the exploration of vibrating structures has leaned on theoretical models followed by experimental testing designed around those models [4]. Vibration analysis of structural systems is essential for understanding the effects of different loading scenarios on the behavior of the structure and identifying parameters such as mass, stiffness, and damping [5]. Recent advancements in data-driven methodologies have presented novel prospects for comprehending the dynamics of vibrating structures directly from measurements [6].

System identification is the process of building an accurate mathematical model of a system using measured data (such as the position over time) and observations [7]. This process is paramount in science and engineering because it provides a means to predict the behavior of a system, even under conditions that are difficult or impossible to replicate in testing. System identification proves critical in comprehending and controlling complex systems [8]; however, systems exhibiting nonlinear and nonstationary vibrations pose significant challenges when identifying their governing dynamics. Conventional identification methods relying on linear assumptions and Fourier-transform-based techniques often fall short [9], resulting in misconstrued dynamical phenomena. Addressing this challenge necessitates new approaches capable of accurately capturing the nonlinear and nonstationary dynamics of such complicated systems [10].

A great collection of methods for nonlinear system identification (NSI) in structural systems are reviewed by [1,2]. The methods commonly employed in NSI can be grouped into three types: (1) parametric methods, (2) non-parametric methods, and (3) semi-parametric methods. Parametric methods focus on identifying values for parameters when the mathematical structure of the model or the physics of the system is already known, for example using, the Non-linear AutoRegressive Moving Average models with eXogenous inputs [11], a Bayesian approach based on the Unscented Kalman Filter method [12], approximate Bayesian computation [13,14], a nonlinear-mode and a polynomial nonlinear state-space system models [15], using time-series optimization [16], a machine learning framework that utilizes moving horizon nonlinear optimization [17], etc. Non-parametric are used when there is no prior knowledge of the underlying physics, e.g., the restoring force surface (RFS) method [18,19], analytical methods, data-driven techniques, and frequency-energy plots [20], neural networks [21], the Nonlinear Identification through eXtended Outputs method with a general high order polynomial for the nonlinear part [22], symbolic regression and genetic programming to identify not only the model but also conservation laws [23]; and leveraging modern computational techniques, a data-driven machine learning approach in combination with symbolic regression [24], infusing physics into a neural network [25], etc. Semi-parametric methods combines the strengths of both parametric and non-parametric methods, [26] introduced a piecewise-linear RFS method, [10] investigated the sparse identification of nonlinear dynamics (SINDy) method [27] in systems that have highly nonlinear elastic or inelastic behavior, [28] replaced in SINDy the sparse regression by a Bayesian framework, [29] used frequency-displacement plots, instantaneous damping curves, and a generic polynomial stiffness model.

Particularly, the SINDy method represents a significant milestone in discovering the governing equations of dynamical systems from measurement data [27]. The process involves employing a dictionary of possible functions that could appear in the system's equations, then a sequentially-thresholded least-squares between the derivative of the collected sample measurements and the dictionary is performed, to identify the fewest terms that represent the data.

SINDy has found application across diverse domains: in fluid mechanics to model low-dimensional complex fluid flows [30]; in heat transfer to identify a reduced model for convection in plasma [31]; in structural engineering to identify nonlinear and hysteric behaviors [10]; in chemical engineering to model and control processes [32]; in biology to represent biological networks [33]; in optics, to identify and correct nonlinear impairments in fiber optic transmissions [34]; and others.

To characterize the nonlinear features of the system, [35] presented an innovative method based on ME, in which a spline-based averaging method was employed to obtain a time-averaged total ME from the averaged kinetic energy (KE) for a linear oscillator with an essentially nonlinear attachment. This total energy is used to define effective measures of stiffness and damping based on the response time series of the linear structure subjected to impulsive loading. These measures implicitly accounted for the effects of the nonlinear attachment on the rate of energy dissipation in the linear oscillator through the variations of the instantaneous energies. The measures also provided insights into the enhancement of the stiffness and damping properties of the system caused by the nonlinear attachment. While effective in the showcased systems, the absence of experimental validation, noisy conditions, and other sources of dissipation, this method may encounter limitations, affecting its accuracy in calculating ME and potentially compromising system identification.

This paper introduces a new data-driven approach for nonlinear SDOF systems, called the energy-based dual-phase dynamics identification (EDDI) method. This method is a semi-parametric method that directly computes the ME from the KE while concurrently deriving the dissipative model of the system. This is possible because the kinetic and mechanical energies are equivalent when the displacement of the oscillator approaches zero. The dissipated energy is then leveraged to compute the ME, which is then used to compute the reformulated Lagrangian. The reformulated Lagrangian is subsequently used to compute the conservative force acting on the oscillator, which is then used to identify a mathematical model for the nonlinear stiffness of the system. As such, the EDDI method identifies both nonlinear damping and stiffness models directly from the measured free response of the oscillator. The method is demonstrated using both simulated and experimentally measured systems with a variety of nonlinearities. Future work will focus on the extension of the proposed method to multi-DOF systems.

## 2. The Proposed Method

Inspired by [35], where they found the link between the ME and the KE, the proposed method is divided into two phases: a damping-model identification phase and a stiffness-model identification phase. The damping-model identification phase focuses on identifying the damping model by computing a model dissipated energy that matches the energy dissipated in the measurements. The second phase uses the obtained dissipated energy to compute the ME of the system, which is then leveraged to obtain the reformulated Lagrangian and the corresponding conservative force acting on the oscillator. A mathematical model for the stiffness nonlinearity is identified by solving a system of linear equations formed by the conservative force and the generalized dispalcement. The following subsections formalizes the exact approach for each phase of the method.

## 2.1. Phase One: Nonlinear Damping Model Identification

This research focuses on the identification of the dynamics of SDOF oscillators experiencing both linear and nonlinear vibrations. The proposed method requires that the user knows the mass of the oscillator and possesses measured free-response data to compute the energies and identify a model for the dynamics of the system. In the case of experimental measurements from accelerometers, the displacements and velocities can be obtained from using numerical integration and filtering as in [36].

The first phase of the proposed method focuses on the calculation of the ME of the system by using the KE, $T(t)$. Considering a generic Duffing oscillator [37] (shown in Fig. 1) with the equation of motion

$$m\ddot{x} + B(x, \dot{x}) + K(x) = F(t), \tag{1}$$

where $B(x, \dot{x})$ represents the internal non-conservative force, $K(x)$ represents the internal conservative force, and $F(t)$ is an externally applied force. In order to propose a generalized method, a Lagrangian framework is used. In this sense, the governing equations are formulated in terms of the generalized coordinate $q$. Hence, the instantaneous ME of the system is

$$E(t) = T(t) + V\big(q(t)\big) = \frac{1}{2}\dot{q}p(t) + V\big(q(t)\big), \; t \geq 0, \tag{2}$$

where, $q(t) = x(t)$, $p(t)$ is the conjugate momentum of $q(t)$, and $V(q(t))$ is the potential energy that depends entirely on the displacement of the oscillator. We assume that the potential energy is zero when the displacement is zero (i.e., $V(0) = 0$), which implies that the displacement, $q(t)$, is measured from static equilibrium. The instantaneous ME can also be expressed as

$$E(t) = E(0) - \int_0^t \dot{q}B(q, \dot{q})d\tau, \; t \geq 0, \tag{3}$$

where the integral term represents the energy dissipated by internal non-conservative force.

Let $\gamma_i, i \in \mathbb{N}$, be the times when the displacement is zero (i.e., $q(\gamma_i) = 0$), such that the potential energy is zero and the ME is equal to the KE

$$E(\gamma_i) = T(\gamma_i) = E(0) - \int_0^{\gamma_i} \dot{q}(\gamma_i)B\big(q(\gamma_i), \dot{q}(\gamma_i)\big)d\tau, \tag{4}$$

such that

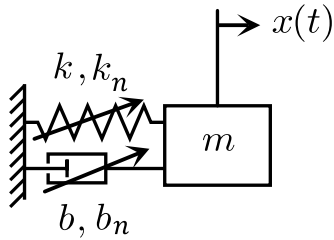

**Fig. 1.** A schematic representation of a Duffing oscillator.

$$E(0) - \int_0^{\gamma_i} \dot{q}(\gamma_i) B\big(q(\gamma_i), \dot{q}(\gamma_i)\big) d\tau - T(\gamma_i) = 0. \tag{5}$$

Note that Eq. 5 represents a discrete-time signal at times $t = \gamma_i$, $i \in \mathbb{N}$ and that the times $\gamma_i$ do not need to be uniformly spaced. Rearranging, we get

$$\int_0^{\gamma_i} \dot{q} B(q, \dot{q}) d\tau = E(0) - T(\gamma_i). \tag{6}$$

Since we may not know $E(0)$, but do know $E(\gamma_0) = T(\gamma_0)$, we shift time to $\gamma_0$, such that $t \geq \gamma_0$ and

$$\int_{\gamma_0}^{\gamma_i} \dot{q} B(q, \dot{q}) d\tau = T(\gamma_0) - T(\gamma_i), \; i \in \mathbb{Z}^+. \tag{7}$$

Note that the right-hand side of Eq. 7 is computed entirely from the measurements while the left-hand side represents the unknown model for the nonlinear damping of the oscillator. As such, the the right-hand side of Eq. 7 provides the necessary information to select a mathematical model for the nonlinear damping. To this end, we propose a model for $B(q, \dot{q})$, such that

$$B(q, \dot{q}) = \sum_{j=1}^{M} b_j B_j(q, \dot{q}), \tag{8}$$

where $b_j$ is a constant coefficient and each $B_j(q, \dot{q})$ contains only a single term. The functions contained in the proposed $B(q, \dot{q})$ may be a combination of polynomial functions or other nonlinear functions, with $M$ denoting the number of proposed terms. Next, we subsitute Eq. 8 into Eq. 7 and expand the result into a set of algebraic equations:

$$b_1 \int_{\gamma_0}^{\gamma_i} \dot{q} B_1(q, \dot{q}) d\tau + \cdots + b_M \int_{\gamma_0}^{\gamma_i} \dot{q} B_M(q, \dot{q}) d\tau = T(\gamma_0) - T(\gamma_i), i \in \mathbb{Z}^+, \tag{9}$$

which can be placed into matrix form and truncated to $i \in [1, N]$ to get

$$\begin{bmatrix} \int_{\gamma_0}^{\gamma_1} \dot{q} B_1(q, \dot{q}) d\tau & \cdots & \int_{\gamma_0}^{\gamma_1} \dot{q} B_M(q, \dot{q}) d\tau \\ \vdots & \ddots & \vdots \\ \int_{\gamma_0}^{\gamma_N} \dot{q} B_1(q, \dot{q}) d\tau & \cdots & \int_{\gamma_0}^{\gamma_N} \dot{q} B_M(q, \dot{q}) d\tau \end{bmatrix} \begin{bmatrix} b_1 \\ \vdots \\ b_M \end{bmatrix} = \begin{bmatrix} T(\gamma_0) - T(\gamma_1) \\ \vdots \\ T(\gamma_0) - T(\gamma_N) \end{bmatrix}. \tag{10}$$

Equation 10 can be written more compactly as

$$Qb = R, \tag{11}$$

and solving this equation produces the model parameters for $B(q, \dot{q})$. Note that the proposed damping model is crucial to identify the dynamics of the model. Finally, from Eq. 3 the estimated ME is obtained as

$$E(t) = E(\gamma_0) - \int_{\gamma_0}^{\gamma_N} \dot{q} B(q, \dot{q}) d\tau, t \in [\gamma_0, \gamma_N]. \tag{12}$$

## 2.2. Phase Two: Nonlinear Stiffness Model Identification

The second phase of the EDDI method focuses on identifying a mathematical model for the internal conservative force and relies on the ME obtained in phase one. After identifying the model for the dissipative force, $B(x, \dot{x})$, the ME is computed by Eq. 12, then the Lagrangian is written as

$$L = T - V = 2T - E. \tag{13}$$

Because $L = L(q(t), \dot{q}(t))$, we arrive at

$$\frac{dL}{dq} = \frac{\partial L}{\partial q}\frac{dq}{dq} + \frac{\partial L}{\partial \dot{q}}\frac{d\dot{q}}{dq}, \tag{14}$$

and

$$\frac{dL}{dq} = \frac{\partial L}{\partial q} + \frac{\partial L}{\partial \dot{q}}\frac{\frac{d\dot{q}}{dt}}{\frac{dq}{dt}}. \tag{15}$$

Then, replacing $\partial L / \partial \dot{q} = p$, we obtain the generalized conservative force as

$$\frac{\partial L}{\partial q} = \frac{dL}{dq} - p\frac{\ddot{q}}{\dot{q}}, \tag{16}$$

To calculate Eq. 16, two considerations are necessary: 1) we take the discrete derivative of the Lagrangian and perform element-wise division to the change in the generalized coordinate, and 2) since we assume that the generalized coordinate and the kinetic energy of the system are known, $p$ can be obtained. The resulting internal conservative force is then plotted as a function of the displacement to form a restoring force plot. The resulting restoring force plot is used to interrogate the dynamics of the system, then a mathematical model is selected based on the observations (e.g., in Cartesian coordinates, $\sum_{n=1}^{5} k_n x^n$). The unknown parameters in the model can then be identified by solving a system of linear equations, using curve-fitting procedures, or other optimization methods to fit the model to the obtained internal conservative force.

## 3. Method Demonstration

### 3.1. Simulated single-degree-of-freedom nonlinear oscillator.

To illustrate the EDDI method, we use simulated data from a Duffing oscillator (DO), adapted from [38], governed by the following equation of motion

$$m\ddot{x} + b\dot{x} + b_{nl}x^2\dot{x} + kx + k_{nl}x^3 = 0, \quad x(0) = 0, \dot{x}(0) = 10, \tag{17}$$

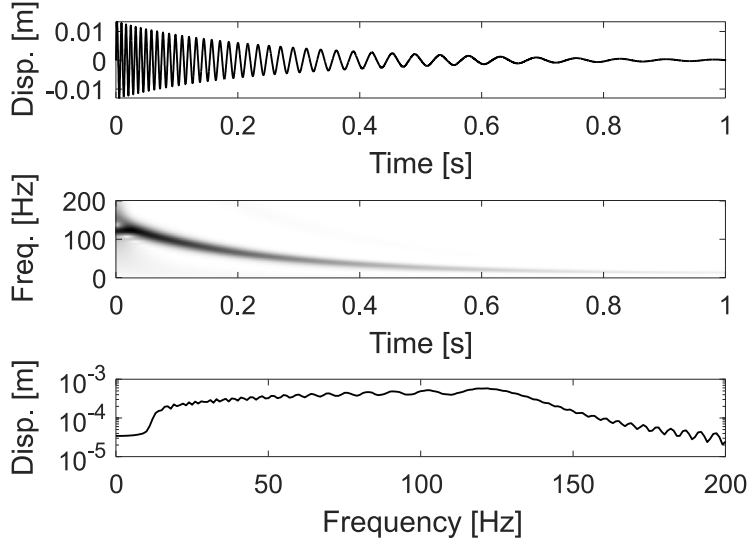

**Fig. 2.** Simulated displacement response of the DO.

where $m = 0.05$ kg, $b = 0.5$ Ns/m, $b_{nl} = 4000$ Ns/m³, $k = 300$ N/m, and $k_{nl} = 3 \times 10^8$ N/m³ are the mass, linear damping, nonlinear damping, linear stiffness, nonlinear stiffness coefficients, respectively, and $x$ is the displacement coordinate. To generate the simulated time series, we use MATLAB® to numerically solve Eq. 17 using the *ode45* solver. The time span is set to $t \in [0,1]$ s, which is long enough for the response to decay to near-zero amplitude, with a sampling rate of $10^4$ Hz. Additionally, the relative and absolute tolerances of the solver were set to $10^{-12}$ and $10^{-16}$, which are more stringent than needed to ensure accurate numerical integration of Eq. 17 without increasing computational costs. Figure 2 shows the simulated time series with the corresponding continuous wavelet transform (CWT) spectrum [39,40] and the Fourier spectrum. The CWT provides a time-frequency representation of the signal content with darker shading representing higher energy content at a particular time and frequency. It is worth nothing that alternative techniques, such as the S-transform [41] or the linear chirplet transform [42] could also be employed to achieve similar depiction. The CWT spectrum in Fig. 2 has been normalized to have a maximum amplitude of 1.

Figure 3(a) shows both the linear and nonlinear damping forces present in the system, which demonstrates that the nonlinear damping is strong for the chosen parameters. Figure 3(b) shows the KE, $T(t)$, the instants when the displacement is zero, $T(\gamma_i)$, and the exact ME. The energy plots reveal a rapid loss of energy in the beginning with a much slower decrease starting around 10 seconds. To capture this behavior, we set the model for the dissipative force as

$$B(x,\dot{x}) = b_1\dot{x} + b_2\dot{x}^2 + b_3\dot{x}^3 + b_4x^2\dot{x}, \tag{18}$$

where the final term represents a mixed displacement-velocity damping term that arises when viscous dampers are excited transversally [43]. By solving Eq. 11, we obtain the parameter values for Eq. 18 and present a comparison of these with the exact values in Table 1. We find that the identified values are close to the theoretical ones. Additionally, we include the values obtained by the SINDy method [27] for the same candidate functions as the EDDI method. We find that SINDy also does a good job of identifying the parameters, but EDDI does a better job in this example.

From these, we can compare and see that the first phase of the EDDI method performs well in identifying the damping model. Note that the assumptions in the model play a key role in the proposed method, such that if the proposed model is incorrectly constructed, the resulting identification will produce a model that does not reproduce the desired behavior. For instance, in our analysis, the identification of parameter values for Eq. 18 heavily relies on the effective solution of Eq. 11, which in turn depends on the appropriate selection and construction of the model library. Integrating Eq. 18, we produce an estimate for the dissipated energy and present a comparison with the exact dissipated energy in Fig. 3(c), where a good agreement is observed. Applying Eq. 12, we estimate the ME for the response and depict a comparison with the exact ME in Fig. 3(d). The results show that the estimated ME agrees well with the exact ME and that the proposed approach can extract the ME directly from free-response signals.

Table 1. Comparison of the coefficients of the dissipative model.

| Coefficient | Exact | EDDI Identified | EDDI Error | SINDy Identified | SINDy Error |
|---|---|---|---|---|---|
| $b_1$ [Ns/m] | 0.5 | 0.49999 | 0.002% | 0.5025 | 0.5% |
| $b_2$ [Ns$^2$/m$^2$] | – | $8.6 \times 10^{-5}$ | – | – | – |
| $b_3$ [Ns$^3$/m$^3$] | – | $-1.2 \times 10^{-7}$ | – | – | – |
| $b_4$ [Ns/m$^3$] | 4000 | 4000.2 | 0.05% | 3923.1 | 1.92% |

We proceed to the second phase by noting that the conjugate momentum for this system is $p = m\dot{x}$. Next, we the compute the conservative force for the response using Eq. 16 and depict the resulting restoring force in Fig. 4. Based on the observed conservative force, we select a polynomial stiffness model using the model $\sum_{n=1}^{5} k_n x^n$, which is reasonable to describe the nonlinear behavior of SDOF oscillator systems [14,44]. Using the chosen model, we identify the unknown parameters by solving the linear system of equations between the proposed candidate of terms and the calculated conservative force. Using the chosen model, we identify the unknown parameters by using the backslash operator (\) in MATLAB®. Nevertheless, other techniques such as curve fitting methods could be used. The exact and identified parameters are provided in Table 2 along with the percent error for each parameter. From the results in Table 2, we find that both $k$ and $k_{nl}$ are correctly obtained, and the proposed method provides a slightly better characterization of the system compared to the model identified by SINDy. Figure 4 shows the comparison of the restoring force for the simulated response and the identified model in the domain $x \in [-0.013, 0.013]$ m. It indeed reveals the nonlinearity due to the coefficient $k_3$. The force calculated by the proposed method is very close to the real one, implying that the proposed stiffness model is correct.

Using the identified model (both damping and stiffness), we simulated the response of the system using the same initial conditions as those used for the exact system. To this end, we present presents a comparison of the displacement time series, CWT spectra, and Fourier spectra for the exact system and identified model in Figure 5. The results show that the identified model accurately reproduces the response of the exact system and that the EDDI method can effectively identify and model the dynamics of the system. We note that, for such a simple system, a strong agreement is expected especially since there is no noise or external influence acting on the simulated system. Nevertheless, the results demonstrate the potential of the EDDI method and its usefulness in identifying mathematical models for strongly nonlinear systems with both nonlinear damping and stiffness.

Table 2. Comparison of the coefficients of the conservative model

| Coefficient | Exact | EDDI Identified | Error | SINDy Identified | Error |
|---|---|---|---|---|---|
| $k_1$ [N/m] | 300 | 317.32 | 5.77% | 361.03 | 20.3% |
| $k_2$ [N/m$^2$] | 0 | $-2.3 \times 10^5$ | – | 581.78 | – |
| $k_3$ [N/m$^3$] | $3 \times 10^8$ | $3.003 \times 10^8$ | 0.1% | $2.987 \times 10^8$ | 0.42% |
| $k_4$ [N/m$^4$] | 0 | $1.08 \times 10^7$ | – | $-5.493 \times 10^6$ | – |
| $k_5$ [N/m$^5$] | 0 | $-4.7 \times 10^{10}$ | – | 0 | – |

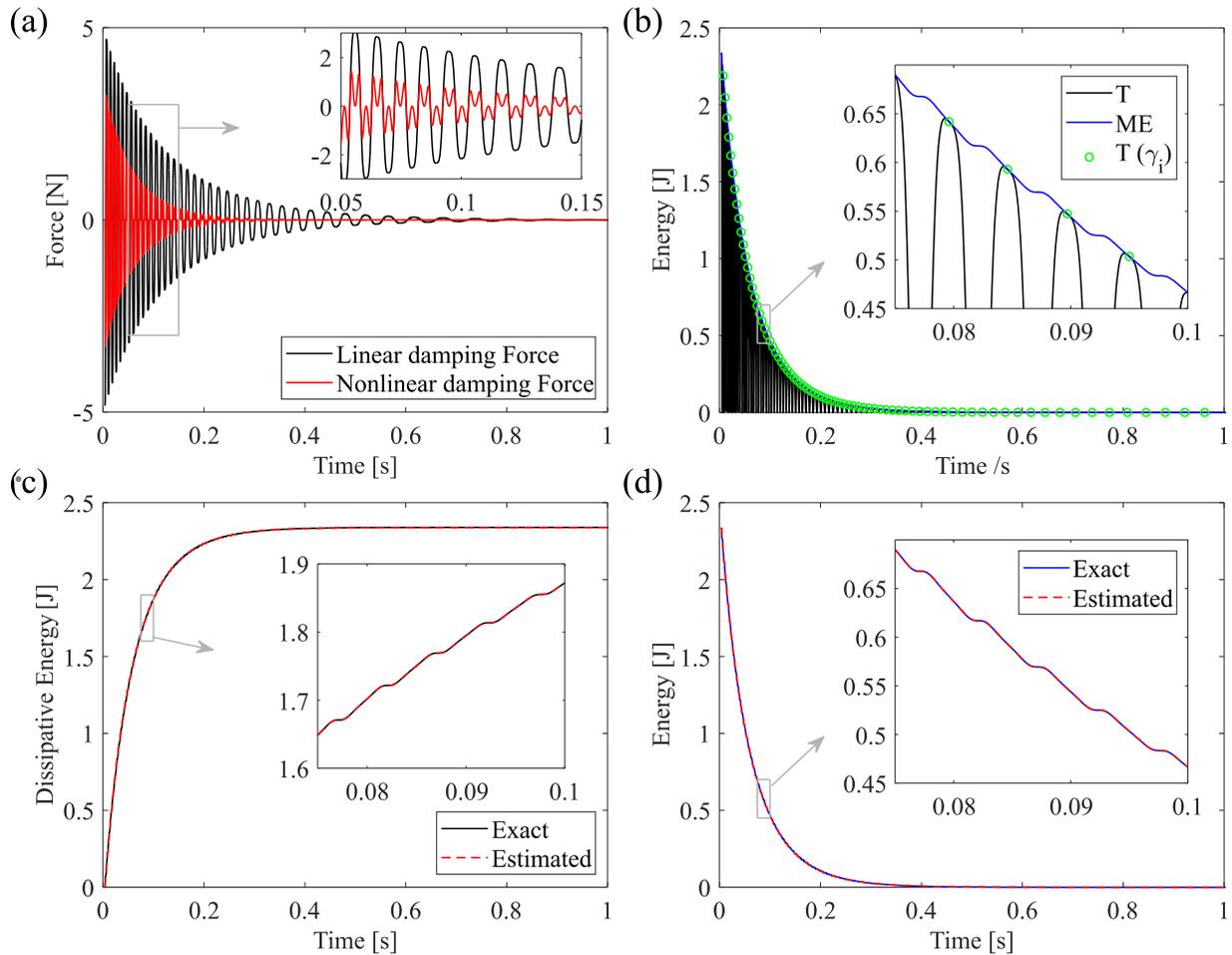

**Fig. 3.** (a) Linear and nonlinear damping forces, (b) Real kinetic and mechanical energies, (c) Comparison of the exact dissipated energy and that estimated by EDDI. (d) Comparison of the exact and estimated mechanical energies.

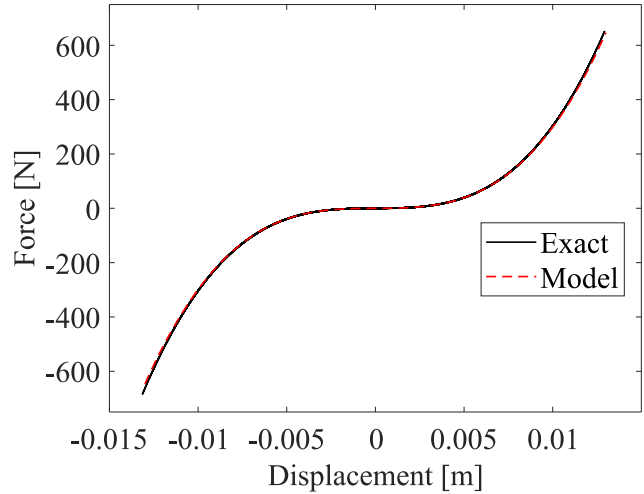

**Fig. 4.** Comparison of the internal conservative forces computed using the exact model and the model identified with the EDDI method shown as restoring force curves.

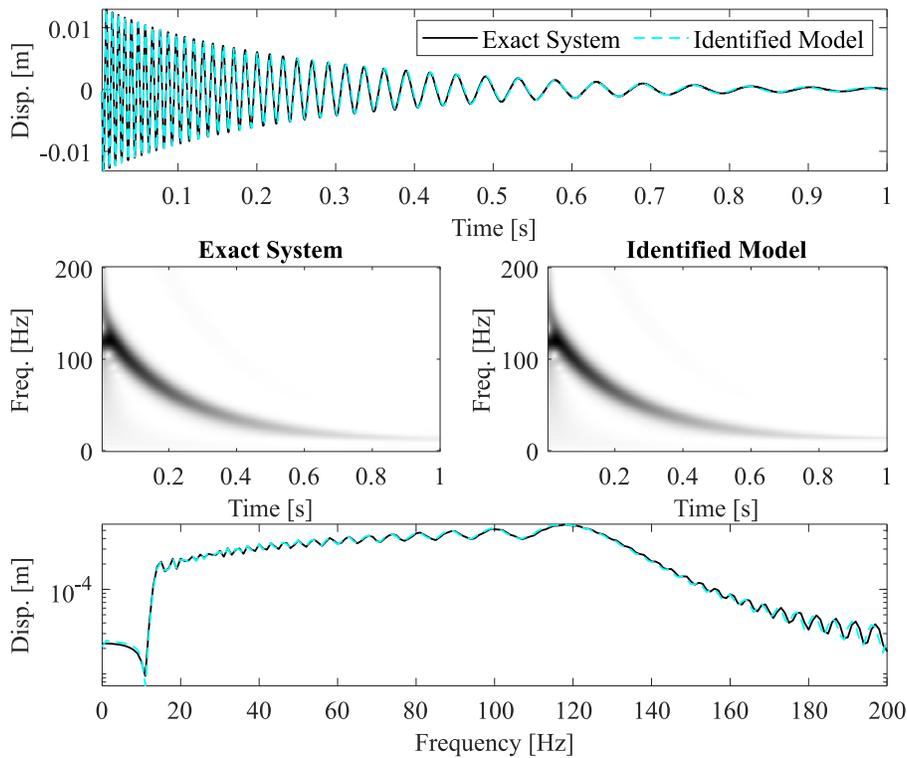

**Fig. 5.** Comparison of the displacement responses, CWT spectra, and Fourier spectra for the exact system and the identified model.

### 3.2. Simulated Damped Pendulum.

To further validate the proposed methodology, we consider the dynamics of a damped pendulum, where the generalized coordinate of motion is the angle $\theta$. The equation of motion of the damped pendulum is

$$ml^2\ddot{\theta} + bl^2\dot{\theta} + mgl\sin(\theta) = 0, \qquad \theta(0) = \pi/2, \dot{\theta}(0) = 0, \tag{19}$$

where $m = 2$ kg, $l = 0.8$ m, $b = 0.1$ kg/s$^2$, and $g = 9.81$ m/s$^2$ are the mass, length, damping coefficient, and acceleration due to gravity, respectively. The response of the pendulum is obtained using *ode45* in MATLAB® for a time span of 100 s with a sampling rate of 100 Hz, and the same tolerances are used as the previous section. Figure 6 depicts the simulated rotation with the corresponding continuous CWT spectrum and the Fourier spectrum. The CWT spectrum reveals a softening-type nonlinearity, which is expected due to the nonlinearity in Eq. 19.

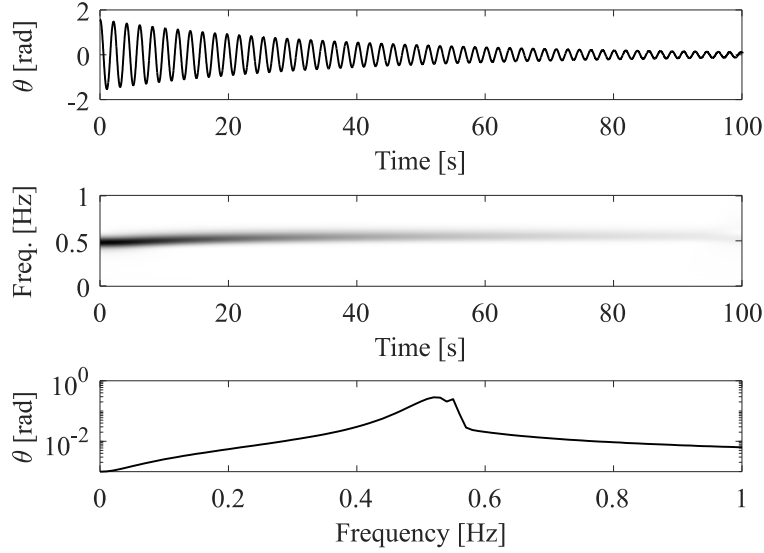

**Fig. 6.** Simulated displacement response of the damped pendulum.

Figure 7(a) shows both the linear damping force of the system, and Fig. 7(b) shows $T(t)$, $T(\gamma_i)$ and the aimed ME. Since this system is a SDOF oscillator, in the first phase, we set the nonlinear damping model as Eq. (18)

$$B(\theta,\dot{\theta}) = b_1\dot{\theta} + b_2\dot{\theta}^2 + b_3\dot{\theta}^3 + b_4\theta^2\dot{\theta}. \tag{20}$$

The values of the solution of Eq. 11 are presented and compared with the theoretical values in Table 3. As we can see there is a well agreement between the calculated and theoretical values. We also present the values obtained by SINDy ($\lambda = 0.005$, and the same candidate functions as used in EDDI), which also provides a close value compared with the theoretical. After integrating Eq. 20, Fig. 7(c) shows the obtained dissipated energy that is similar to the real dissipated energy. Then, by applying Eq. 12 we can compute the approximate ME, $\tilde{E}(t)$, as shown in Fig. 7(d), from which the $\tilde{E}(t)$ resembles the real ME.

For the second phase, we first set $p_\theta = ml^2\dot{\theta}$ [38], then we obtain the conservative force using Eq. 16. After that, we use the stiffness model: $\sum_{n=1}^{5} k_n\theta^n$, and finally, for this model we solve the linear system of equations between the obtained conservative force and the angular displacement. The theoretical values corresponding to the Taylor expansion of the sine function and obtained values are shown in Table 4, from which the values are well approximated. In this case, the SINDy method provides better accuracy compared with the proposed method. Figure 8(a) depicts the loops for both the simulated and the obtained restoring forces in the domain $\theta \in$ [−1.5,1.5] rad, which show relatively good agreement for the nonlinear sine term. Figure 8(b) presents the solution of Eq. 19 for both the real and the obtained by EDDI, with the corresponding CWT and Fourier spectra, in which we can see that the proposed method can also capture the dynamics of this system that is best analyzed in polar coordinates. Indeed, there is a good agreement due to the simulated signal of Eq. 19 was used to find the parameters of the system.

Table 3. Comparison of the coefficients of the dissipative model

| | | EDDI | | SINDy | |
|---|---|---|---|---|---|
| Coefficient | Exact | Identified | Error | Identified | Error |
| $b_1 = bl^2$ [Nms/rad] | 0.064 | 0.064006 | $9.4 \times 10^{-3}$% | 0.063984 | 0.025% |
| $b_2$ [Nms$^2$/rad$^2$] | 0 | $-9.4 \times 10^{-5}$ | – | 0 | – |
| $b_3$ [Nms$^3$/rad$^3$] | 0 | $-3.2 \times 10^{-5}$ | – | 0 | – |
| $b_4$ [Nms/rad$^3$] | 0 | $1 \times 10^{-3}$ | – | 0 | – |

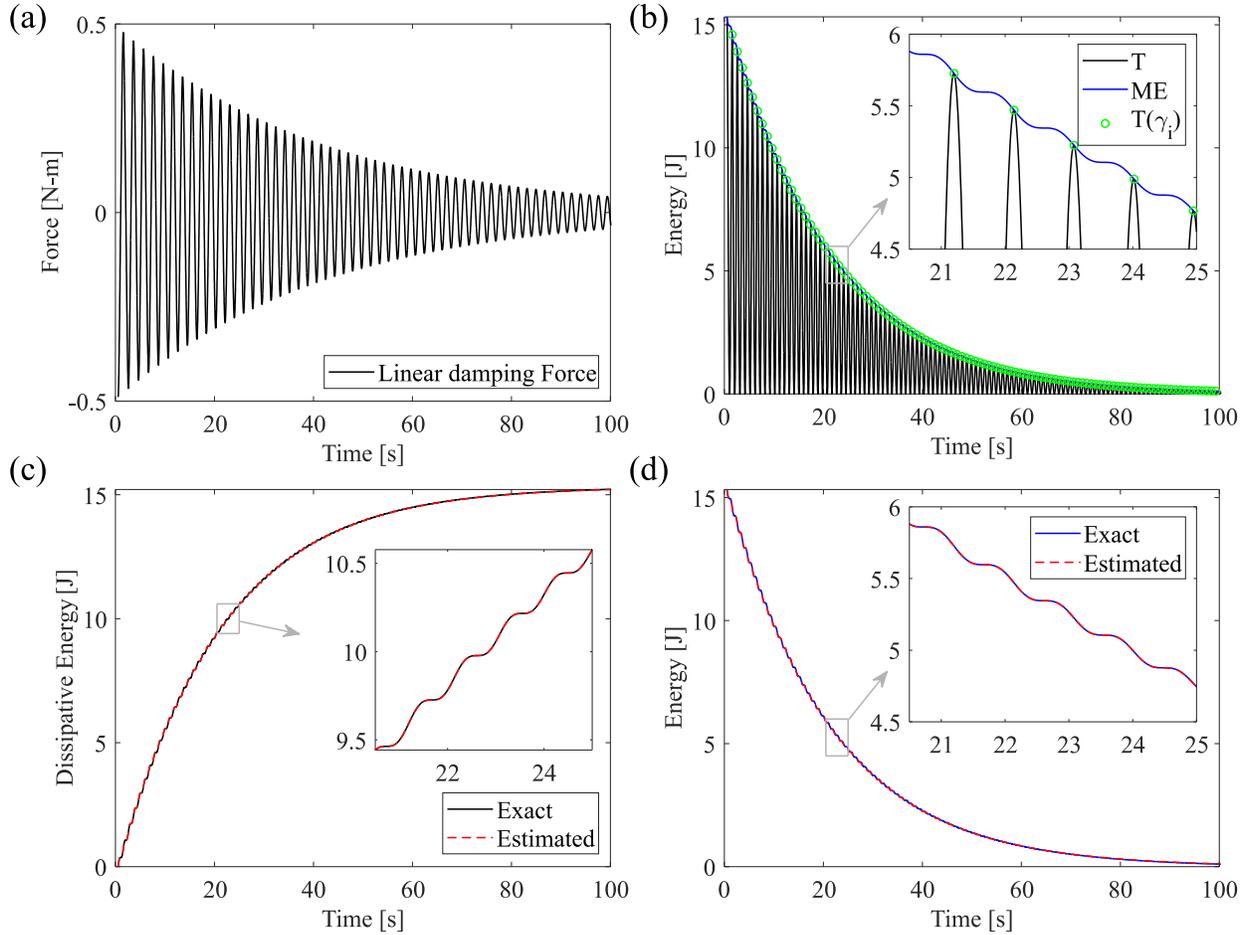

**Fig. 7.** (a) Linear and nonlinear damping forces, (b) Real $T$ and ME, Comparisons of: (c) Real dissipated energy computed by EDDI, (d) Real and calculated ME.

Table 4. Comparison of the coefficients of the conservative model

| | | EDDI | | SINDy | |
|---|---|---|---|---|---|
| Coefficient | Exact | Identified | Error | Identified | Error |
| $k_1 = mgl$ [N/rad] | 15.696 | 15.652 | 0.28% | 15.691 | 0.032% |
| $k_2$ [N/rad$^2$] | 0 | 0.43 | – | 0 | – |
| $k_3 = mgl/3!$ [N/rad$^3$] | −2.616 | −2.58 | 1.34% | −2.605 | 0.42% |
| $k_4$ [N/rad$^3$] | 0 | −0.167 | – | 0 | – |
| $k_5 = mgl/5!$ [N/rad$^5$] | 0.1308 | 0.112 | 14.46% | 0.12 | 7.83% |

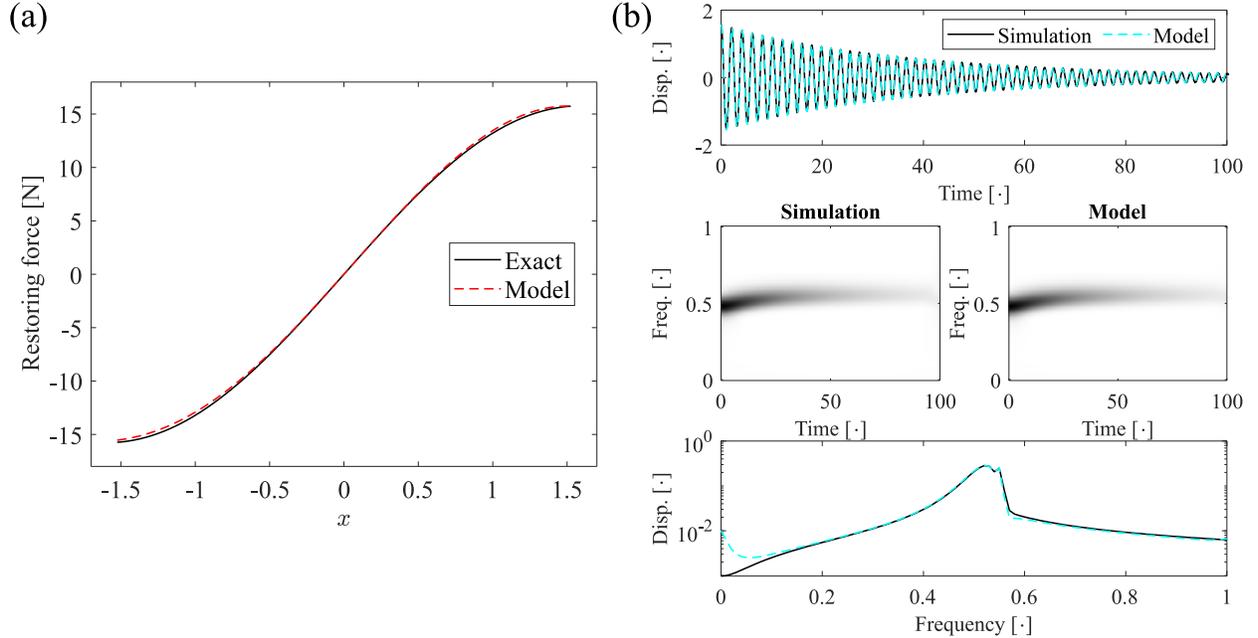

**Fig. 8.** (a) Internal conservative forces computed using the EDDI method and the exact forces, (b) displacement responses, CWT spectra, and Fourier spectra for the exact system and the identified model.

### 3.3. Experimental Demonstration

### 3.3.1 Duffing Oscillator with Smooth Nonlinearity

To demonstrate the EDDI method in real vibration signals, we used an experimental system that resembles a DO with thin flexures and thin wires as springs. The DO (Fig. 9) is grounded to an optical table using an aluminum plate, two rectangular beams, hollow steel flexures and 1/4"-20 UNC bolts. The steel flexures, which serve as linear springs, are bolted to the rectangular beams and the mass. The length, width, and thickness of the steel flexures are 0.1524 m, 0.114 m, and 0.0003 m, respectively. The mass is constructed using a 0.1524 m × 0.1524 m × 0.0127 m aluminum plate and has a mass of 0.815 kg. In addition to the linear springs, the wires introduce both linear and strongly nonlinear stiffnesses that result from a combination of pre-tension and bending. The two parallel steel wires each have a diameter of 0.00074 m and are clamped to C-shaped aluminum brackets at each end and to the center of the DO mass plate using 10-32 UNF set screws.

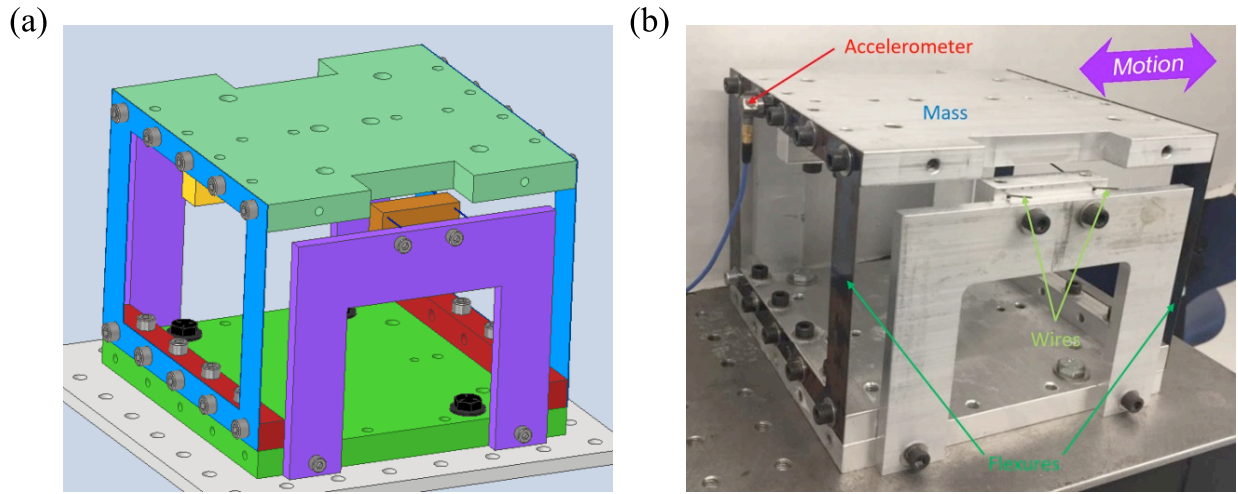

**Fig. 9**. (a) CAD model, (b) Instrumented Duffing oscillator with flexures, wires, mass, C-shaped brackets, and the direction of motion for the nonlinear oscillator.

We performed a series of experimental measurements by applying an impact from a PCB Piezotronics modal impact hammer (model 086C03) to the DO, and we measured the resulting transient acceleration using a PCB accelerometer (model 353B15) with a nominal sensitivity of 1 mV/(m/s$^2$). The measurements were performed using a sampling rate of 19200 Hz using HBM MX1601B hardware and catman data acquisition software. The response was measured for a total duration of 45 s with a pre-trigger length of 1 s. From the measurements, the velocities of the DO are computed by numerically integrating the accelerations and high-pass filtering the output using a third-order Butterworth filter with a cutoff frequency of 2 Hz. The displacement responses were computed applying the same procedure to the corresponding velocities.

Figure 10(a) shows the force applied to the system whereas Fig. 10(b) presents the displacement response, the CWT spectrum, and the Fourier spectrum. The resulting response is nonlinear, reflected on how the frequency varies in time, indicating that the nonlinearity participates strongly in the response for the first 10 s. Figure 10(c) depicts the KE and the corresponding $T(\gamma_i)$. For this system, we propose the following dissipative model for the first phase

$$B_{DO}(t) = b_1 \dot{x} + b_2 x^2 \dot{x} + b_3 \dot{x}^3. \qquad (21)$$

By solving Eq. 11, we extract the following values: $b_1 = 0.196$ Ns/m, $b_2 = -20441$ Ns/m$^3$, and $b_3 = 9.96$ Ns$^3$/m$^3$. It is essential to note that the choice of models extends beyond what we propose here. Interested readers have the freedom to explore and experiment with various other models based on their knowledge, creativity, and intuition. By integrating Eq. 21, we obtain the estimated dissipated energy as shown in Fig. 10(d). Around 5.5 s, we observe a minor discrepancy between the calculated dissipated energy and the experimental value. This discrepancy is attributed to noise in the signal, which adds complexity to the analysis. Nevertheless, the estimated ME calculated by Eq. 12, depicted in Fig. 10(e), appear reasonable given compared to the KE.

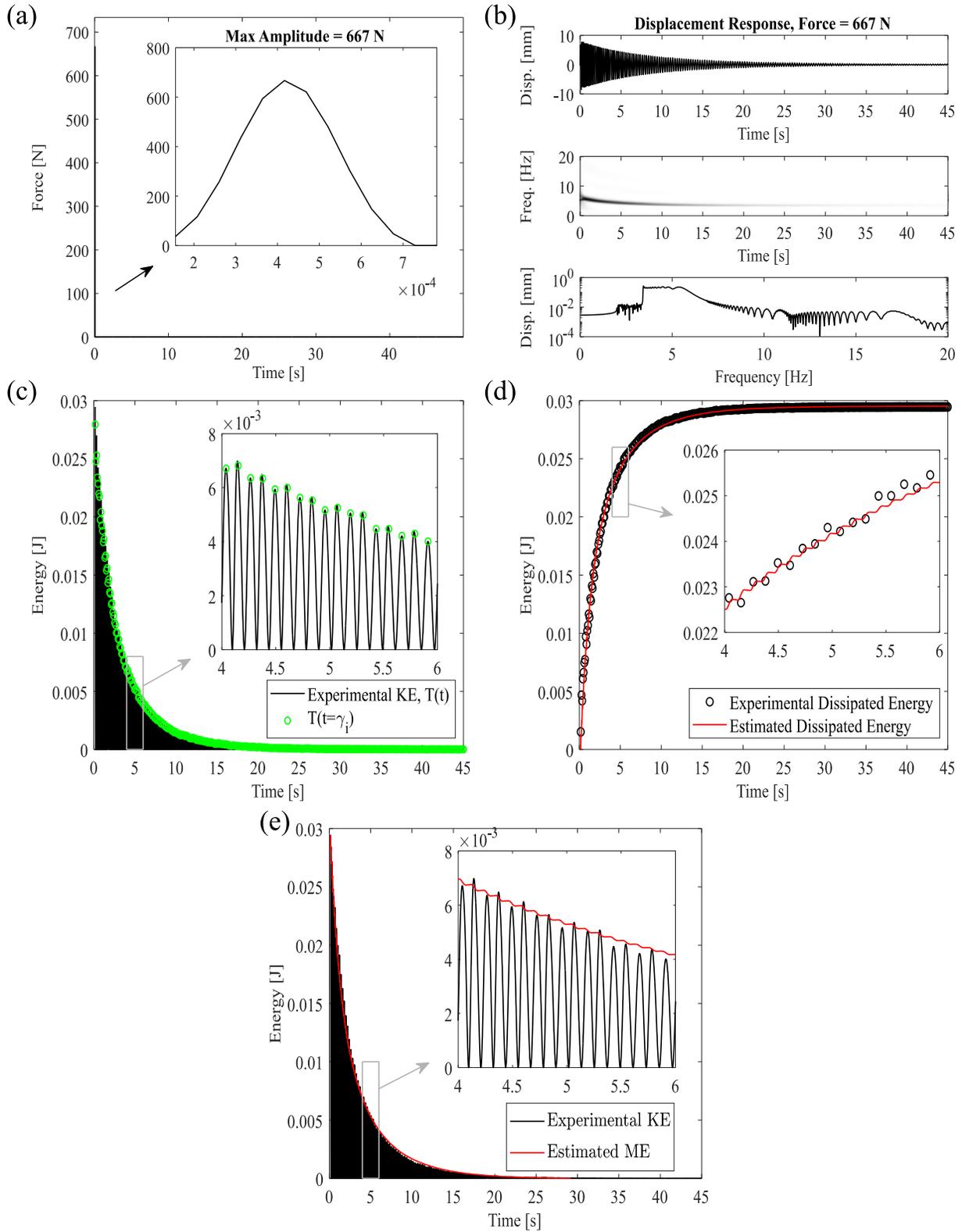

**Fig. 10.** (a) The experimentally measured excitation; (b) the displacement response of the DO; (c) the experimental KE and $T(\gamma_i)$; (d) Experimental and estimated dissipated energies; and (e) comparison of the experimental KE and the estimated ME.

For the second phase, we first compute the Langrangian using the KE and estimated ME using Eq. 13. Second, we set $p_x = m\dot{x}$ and obtain the conservative force using Eq. 16

$$K_{DO}(x) = \sum_{n=1}^{5} k_n x^n. \tag{22}$$

We identify the underlying parameter solving the linear system of equations between the model in Eq. 22 and the displacement. We obtain the following coefficients: $k_1 = 390$ N/m, $k_2 = -1.4 \times 10^3$ N/m², $k_3 = 1.3 \times 10^7$ N/m³, $k_4 = -5.5 \times 10^7$ N/m⁴, and $k_5 = 6.2 \times 10^{10}$ N/m⁵. Figure 11(a) shows both the experimental conservative force (obtained in Eq. 16) and the identified model (obtained in Eq. 22) in the time series. Notably, a good match is observed between the two. Figure 11(b) displays the restoring-force plot, revealing the nonlinear nature of the conservative forces.

Finally, the models for the dissipative and conservative forces are combined to produce the equation of motion of the system as

$$m\ddot{x} + b_1\dot{x} + b_2 x^2 \dot{x} + b_3 \dot{x}^3 + \sum_{n=1}^{5} k_n x^n = F(t), \tag{23}$$

where $F(t)$ is the externally applied force (i.e., the impact from the modal hammer). To validate the model given in Eq. 23, we numerically integrate it using *ode45* in MATLAB® with zero initial conditions and the measured force applied up to a time of 0.01 s. Additionally, the *interp1* function in MATLAB® is used to interpolate times between the sampled points in the measured force signal. Figure 11(c) compares the response of the experimental system and the identified model. The results demonstrate a good agreement between the measured and simulated responses in the time series, CWT spectra, and Fourier spectra.

In addition to the EDDI method, we also applied the SINDy method with the threshold set to 0.05 and the same damping and stiffness models as used in the EDDI method. SINDy resulted in the following parameters for the dissipative model: $b_1 = 0.39$ Ns/m, $b_2 = -1.2 \times 10^5$ Ns/m³, and $b_3 = 37.6$ Ns³/m³. For the stiffness model, SINDy produced the following values: $k_1 = 364.8$ N/m, $k_2 = 329.7$ N/m², $k_3 = 1.6 \times 10^7$ N/m³, $k_4 = -8.7 \times 10^7$ N/m⁴, and $k_5 = 0$ N/m⁵. Figure 11(d) presents a comparison of the displacement time series, CWT spectra, and Fourier spectra for the measured system and that identified using SINDy. As can be seen, SINDy underestimates the linear stiffness and overestimates the damping, such that the simulated response does not match the measured response.

A more robust validation of the EDDI method is obtained by using the identified model to reproduce measurements that were not used in the identification. To this end, we present comparisons for impacts of 246 N and 1110 N in Figs. 12(a) and (b), respectively. For an impact of 246 N, the resulting response is strongly nonlinear, but weaker than the measurement case used for the identification. There is a good agreement between the model and experiment based on a comparison of the time series. The amplitude of the simulated response is slightly smaller than that of the measured response, but the overall trends are well reproduced. This indicates that the damping model accurately reproduces the dissipation at this amplitude whereas the stiffness model may deviate slightly. Similar results are observed for the comparison of the CWT and Fourier spectra. For an impact of 1110 N, shown in Fig. 12(b), the response exhibits a substantially larger change in frequency content as observed in the CWT and Fourier spectra. In terms of the time-series, the simulated response reproduces the experimental measurements except for a slight difference in amplitude that is greatest around a time of 3 s. The Fourier spectra reveals that the nonlinearity identified in the model is not as strong as that in the experiment; however, the overall agreement is good.

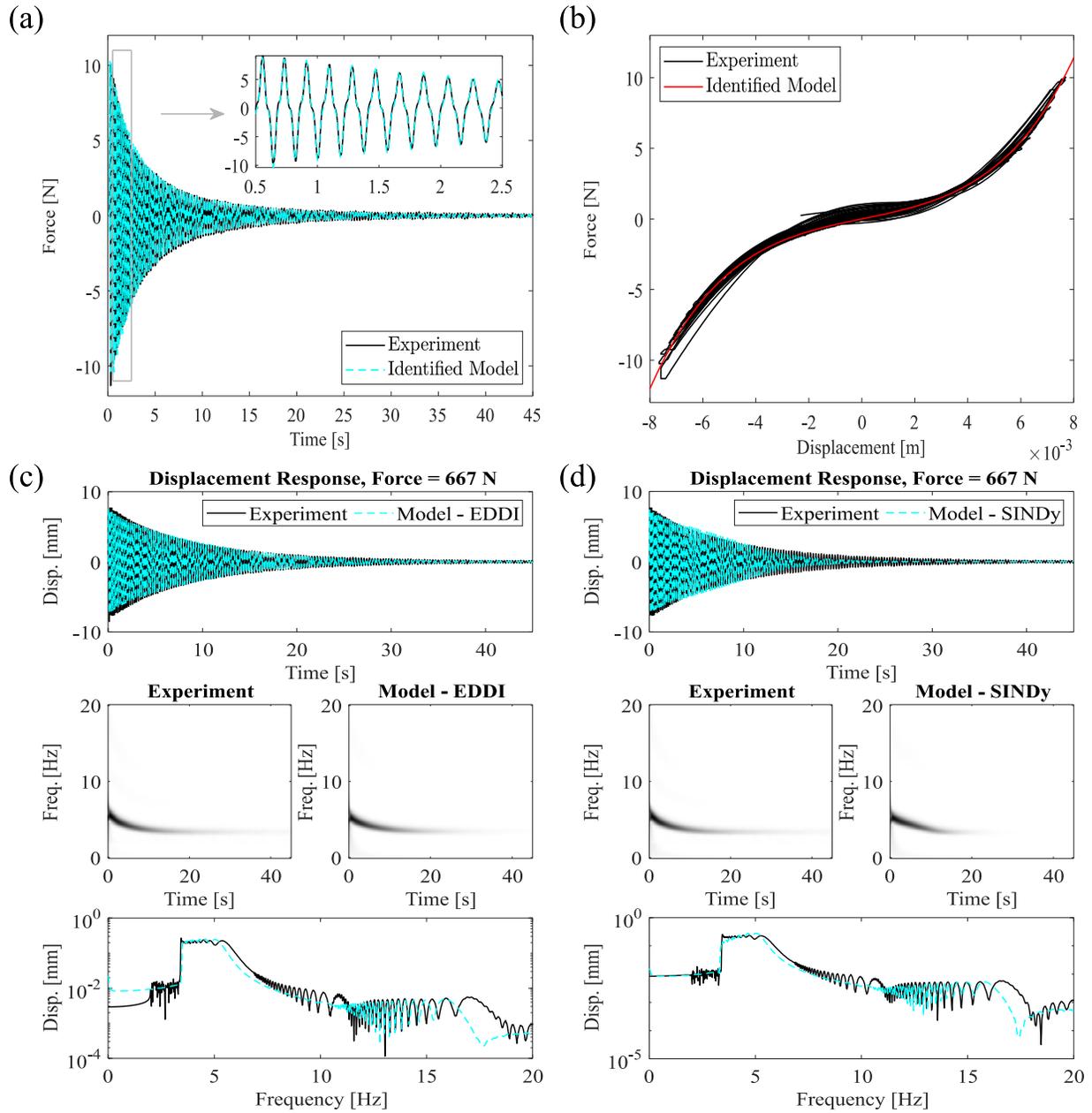

**Fig. 11.** Internal conservative forces computed using EDDI method and the experimental force: (a) Time series, (b) Force-displacement; (c) Validation of the proposed method, (d) result obtained by SINDy method.

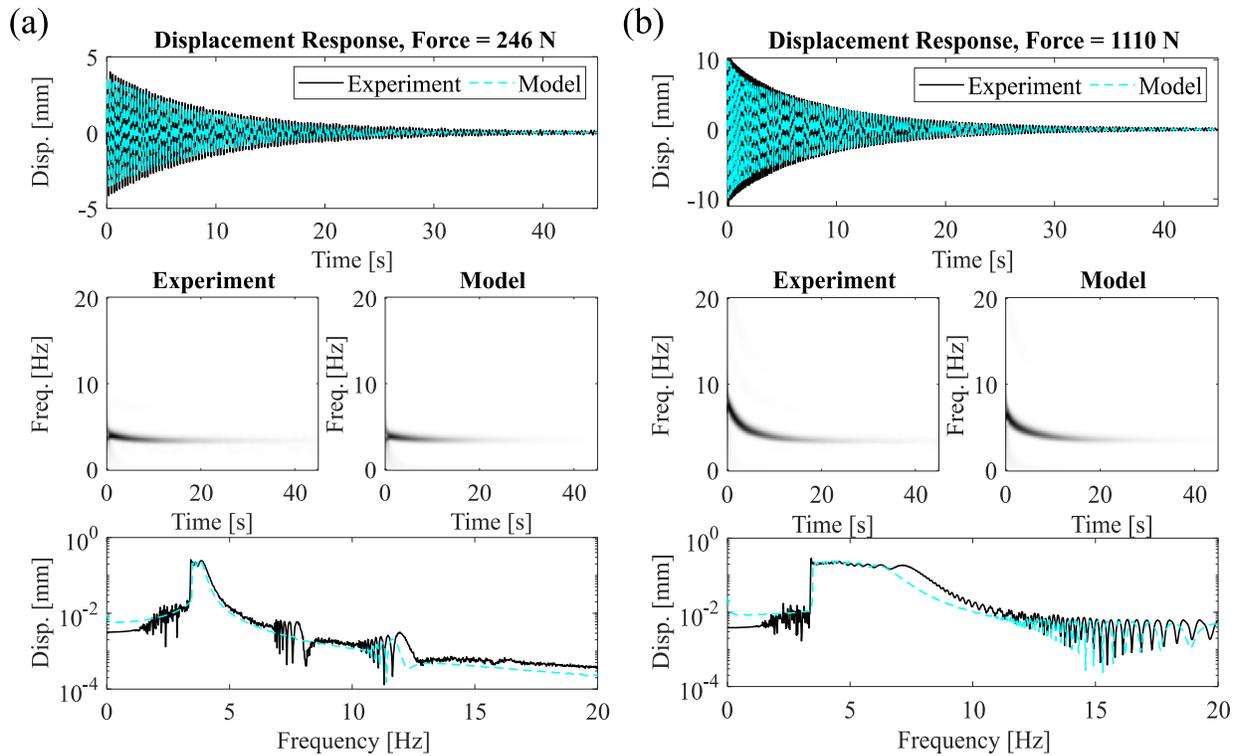

**Fig. 12.** Experimentally measured excitation and displacement responses of the DO for impact amplitudes of (a) 246 N and (b) 1110 N.

### 3.3.2. Duffing Oscillator with Softening-Hardening Nonlinearity

We demonstrate the EDDI method in a second experiment, which consists of a DO with elastic struts that introduce strong damping and stiffness nonlinearity (see Fig. 13). The mass of the DO is suspended by steel flexures using 6-32 UNC screws. The mass is a cuboid made from steel with dimensions of 0.0305 m × 0.0305 m × 0.0127 m and a mass of 0.088 kg. Each end of the steel flexures is grounded to a rectangular, steel beam, which is bolted to an aluminum plate. This plate is grounded with 10-32 UNC screws to L-brackets that are bolted to an optical table. The active length and width of the steel flexures are 0.0702 m and 0.0103 m, respectively. The thickness of each steel flexure is 0.000127 m. The setup of the experimental measurements is the same as the previous subsection, except that a different accelerometer (model 352C03) is used. This accelerometer has nominal sensitivities of 1 mV/(m/s$^2$) and less noise at low frequencies compared to the accelerometer used in the previous subsection. The elastic struts, made of urethane rubber, have a width of 0.0221 m, a height of 0.0193 m, and a thickness of 0.00483 m. A detailed explanation of them is found in [45]. The transient responses were measured for 4 seconds, and the accelerations were integrated then filtered to obtain the velocities and displacements just as in the previous subsection. In this case, the cutoff frequency was set to 3 Hz.

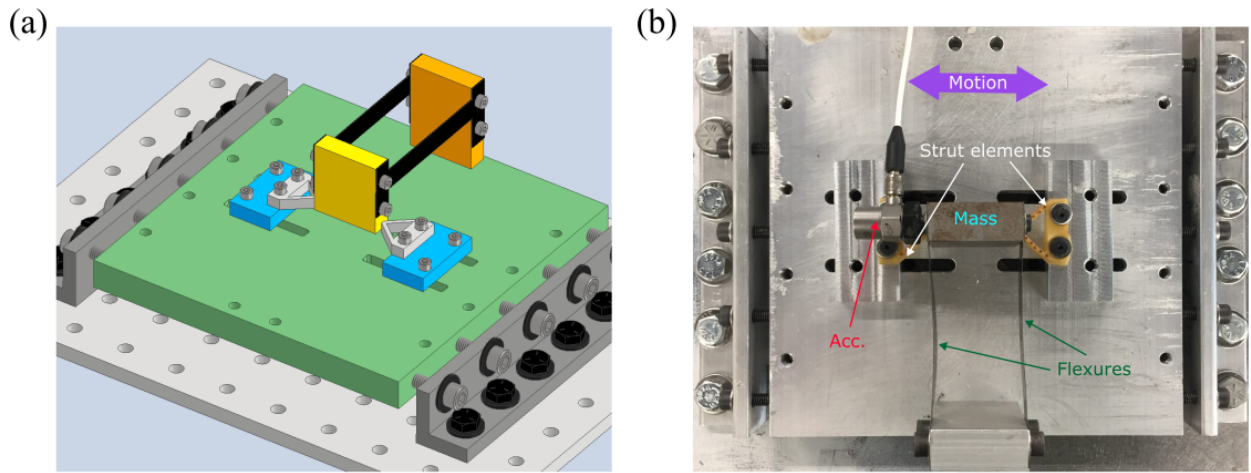

**Fig. 13.** (a) CAD model, (b) Top view of the instrumented Duffing oscillator with flexures, strut elements, mass, and the direction of motion for the nonlinear oscillator.

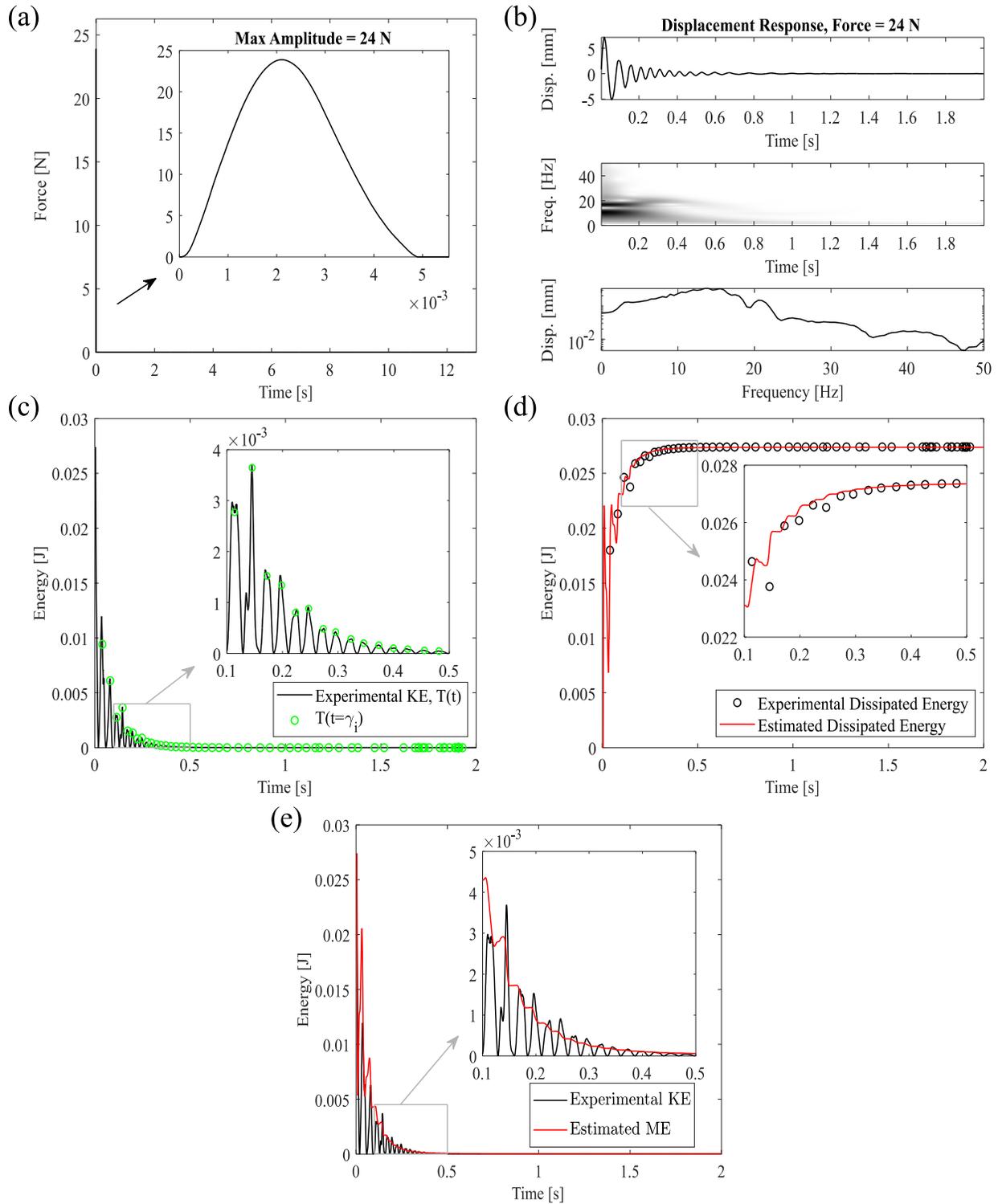

**Fig. 14.** (a) Experimentally measured excitation, (b) Displacement response of the DO-ES, (c) experimental KE and $T(\gamma_i)$, (d) Experimental and modeled dissipative energies, (e) Calculated $ME$, (f) Obtained $V$.

Figure 14(a) depicts the external excitation applied to the system, while Fig. 14(b) presents the corresponding displacement response, the CWT spectra, and the Fourier spectrum. The resulting response is strongly nonlinear, a result of using strut elements, which implies that the response dies out in about 1 s. Figure 14(c) shows both $T$ and $T(\gamma_i)$. For the first phase, we use the same dissipative model as the previous section as

$$B_{ES}(t) = b_1 \dot{x} + b_2 x^2 \dot{x} + b_3 \dot{x}^3 \qquad (24)$$

By solving Eq. 11, we get the values: $b_1 = 0.48$ [Ns/m], $b_2 = -4.07 \times 10^5$ Ns/m$^3$, and $b_3 = 27.7$ Ns$^3$/m$^3$. The estimated dissipated energy is depicted in Fig. 14(d). As in the previous case, due to the presence of noise, the calculated dissipated energy deviates slightly from the experimental values. The ME is estimated using Eq. 12 and the results is depicted in Fig. 14(e). One can notice that in this case, the ME is harder to obtain not only because the elastic strut elements present a strong nonlinearity and the system noise, but also there is a clearance effect when the mass is not in contact with them.

For the second phase, first, the Lagrangian is calculated by Eq. 13. Then, by utilizing $p_x = m\dot{x}$, the conservative force is computed by Eq. 16. To the result of the numerical derivative, the smoothing operation is applied with a window size of 100 data. After that, the stiffness model used is

$$K_{ES}(x) = \sum_{n=1}^{5} k_n x^n. \qquad (25)$$

For the model in Eq. 25, we solve the linear system of equations between the estimated conservative force and the displacement resulting in the following values: $k_1 = 1064$ N/m, $k_2 = 9257$ N/m$^2$, $k_3 = -4.4 \times 10^7$ N/m$^3$, $k_4 = -9.4 \times 10^8$ N/m$^4$, and $k_5 = 8.49 \times 10^{11}$ N/m$^5$. Using SINDy, we obtain the following parameters for dissipative model: $b_1 = 0.75$ Ns/m, $b_2 = 1.98 \times 10^4$ Ns/m$^3$, and $b_3 = 5.99$ Ns$^3$/m$^3$. For the stiffness model, we get $k_1 = 994.1$ N/m, $k_2 = -1.92 \times 10^4$ N/m$^2$, $k_3 = -1.8 \times 10^7$ N/m$^3$, $k_4 = 1.9 \times 10^9$ N/m$^4$, and $k_5 = 0$ N/m$^5$. Figure 15(a) shows both the experimental and model conservative forces as a time series and an acceptable agreement between the two is observed. The restoring-force plots are provided in Fig. 15(b) and these display the common behavior of this system, softening and hardening phenomenon, showing that we capture the overall dynamics of the system. Then, the resulting governing equation of motion is

$$m\ddot{x} + b_1\dot{x} + b_2 x^2 \dot{x} + b_3 \dot{x}^3 + k_1 x + k_3 x^3 + k_5 x^5 = F(t). \tag{26}$$

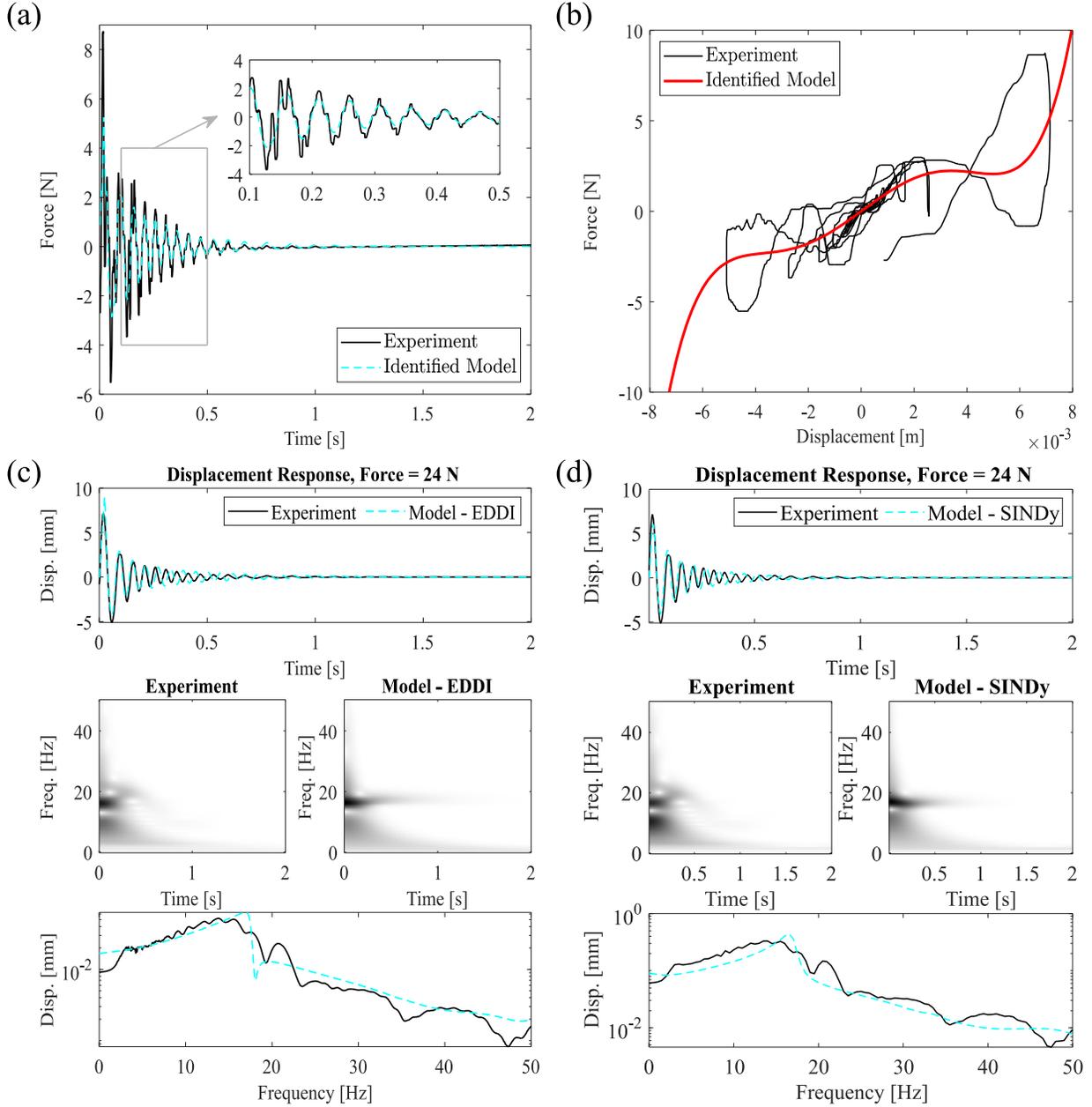

**Fig. 15.** Comparison of the internal conservative forces computed using the EDDI method (experiment) and the identified model shown as (a) time series and (b) restoring-force plot. Comparison of the measured and simulated responses for (c) the model identified using EDDI and (d) the model identified by SINDy.

Figure 15(c) shows the validation of Eq. 26, where there is a fair match between the measured time series and the solution of Eq. 26. This indicates that the identified parameters accurately capture the dynamics of the measured system. When using the SINDy method (see Fig. 15(d)), the resulting model does not accurately capture the softening-hardening behavior as observed in the CWT spectra.

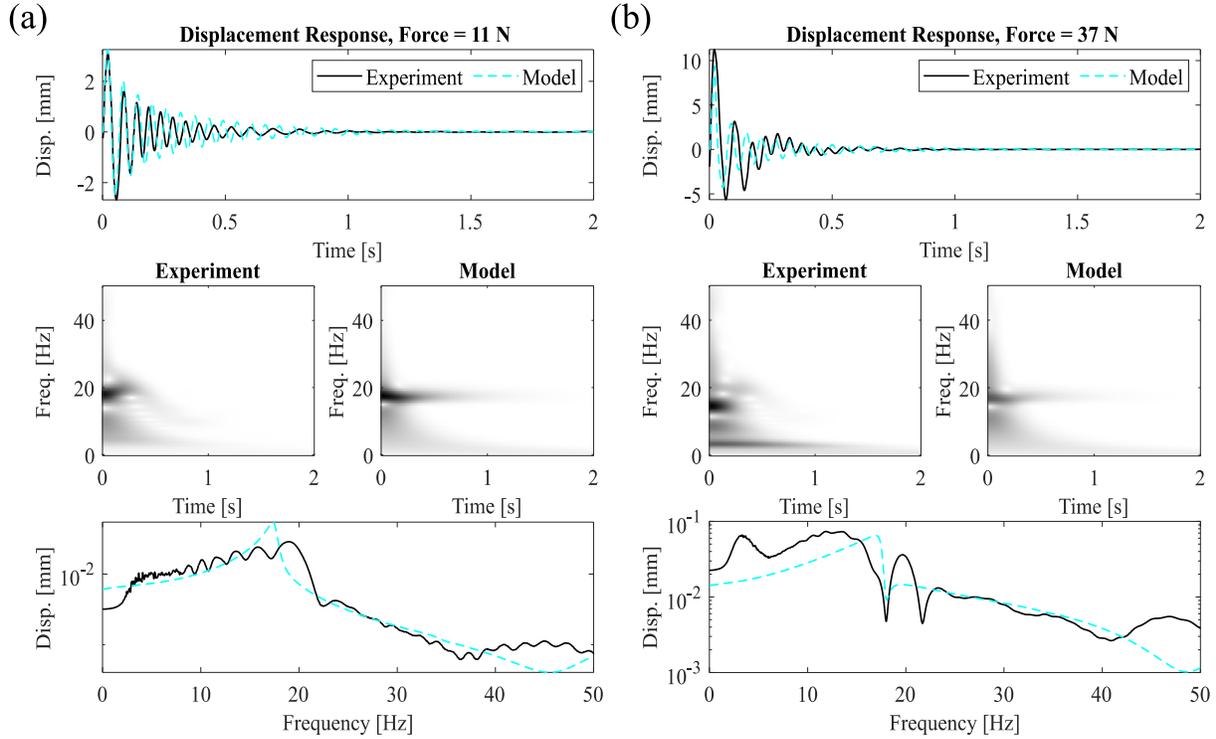

**Fig. 16.** Experimentally measured excitation and displacement responses of the DO-ES for impact amplitudes of (a) 11 N and (b) 37 N.

Figure 16(a) displays the validation for two measurement cases not used in the identification. For an impact of 11 N, the resulting response shows that the frequency rapidly decreases in time from 20 to 10 Hz. For this case, at the beginning our proposed method captures the behavior well. When the impact is 37 N (Fig. 16(b)), the resulting response does not decay exponentially. As a result, the identified model does not agree with the measured response as well as it does for lower amplitude impacts.

## 4. Discussion and Open questions

The results presented in this work demonstrate the power of the novel EDDI method for identifying governing equations of motion for SDOF systems using only measured mass properties and transient response data. The method also relies on basic knowledge that allows one to know the momentum of the system and to interpret the estimated dissipated energy and the estimated conservative force, such that mathematical models can be proposed for each of them. The EDDI technique is applicable to systems that possess smooth nonlinearities. The method is intended for use with free-response data such that large frequency transitions are captured in a single measurement.

When dealing with noisy data, the numerical differentiation performed to compute Eq. 16 is inherently sensitive to noise. Thus, small fluctuations in the data can lead to significant variations in the estimates of the conservative forces (i.e., local infinities arise when the velocity becomes zero). In that regard, careful consideration of smoothing, and appropriate differentiation methods [46] is essential to obtain the conservative force. Additionally, all the systems investigated in this work consisted of smooth nonlinearities and the applicability to non-smooth systems remains to be investigated. The authors suspect that the method could apply directly to clearance-type nonlinearities but would need to be extended to capture systems with rigid-body impacts. This will be explored in future work.

While the EDDI procedure is powerful when dealing with SDOF systems, the current framework cannot be applied to multi-DOF (MDOF) systems. One potential approach to extend the EDDI method to MDOF systems is to apply signal decomposition algorithms to extract individual modes, then apply the EDDI method to the individual modes and combine the identified models. The authors plan to investigate this approach in future work, but also welcome other contributions and ideas from the community at large.

## 5. Concluding Remarks

This research presented a new, data-driven method for identifying the dynamics of SDOF systems called the energy-based dual-phase dynamics identification (EDDI) method. The EDDI method relied entirely on the mass and measured transient response of a SDOF Duffing oscillator to derive and identify a mathematical model for both the underlying stiffness and damping. The EDDI technique was applied to both simulated and experimentally measured responses of Duffing oscillators. In the first phase of the EDDI technique, the measured responses were used to obtain the KE, which was then used to estimate the dissipated energy. Next, a damping model was proposed and the parameters were obtained by fitting the model to the estimated dissipated energy. The resulting model allows one to then estimate the ME of the measured response. In phase two, the Lagrangian is calculated from the kinetic and mechanical energies. The Lagrangian is then used to compute the conservative force. Based on the experimental conservative force, the analyst proposes a stiffness model and identifies the parameters by solving a system of linear equations. For the simulated systems, the identified models were validated through comparisons between the identified and the ability of the model to reproduce the simulated response. For the experimental systems, the method was validated by using the identified model to simulate the response of the system for excitations not used for the identification. The results demonstrated that the EDDI method can produce mathematical models that accurately capture the dynamics of SDOF systems with strong, smooth nonlinearities.

## Acknowledgments

CL is grateful for the support of the Fulbright Program through his Fulbright scholarship. This research was supported by the Air Force Office of Scientific Research Young Investigator Program under grant number FA9550-22-1-0295.